\documentclass [12pt]{article}
\usepackage {graphicx}
\usepackage {longtable}

\begin{document}
\begin{center}
\textbf{On one inverse problem relatively domain for the plate}
\end{center}

\begin{center}
\textbf{Gasimov Y.S.}
\end{center}

\begin{center}
Institute of Applied mathematics Baku State University, Z. Khalilov 23,
\end{center}

\begin{center}
AZ1148 Baku, Azerbaijan, e-mail: .\underline {ysfgasimov@yahoo.com}
\end{center}

\bigskip

\textbf{Abstract.} Inverse problem relatively domain for the plate under
across vibrations is considered. The definition of s-functions is
interoduced. The construction for defining of the domain of the plate by
given s-functions is offered.

\bigskip

Plates are elements of the constructions, which are widely used in various
technical solutions. In this connection investigation of the different
characteristics of the plates is one of the actual problems of the optimal
projecting theory [1]. We consider the problem of finding of the form of the
plate under across vibrations by given characteristics of the system-
inverse problem relatively domain.

In traditional inverse spectral problems by given some experimental data
(scattering data, normalizing numbers etc.) the potential is reconstructed
or the necessary and sufficient conditions are proved providing unequivocal
determination of the seek functions [2].

In differ from this, inverse spectral problem relatively domain has another
specification. Firstly, these problems require to find rather a function but
domain. Secondly, the choice of data (results of the observation),
sufficient for reconstruction of the domain is also enough difficult problem
[3-4].

Let $D$ be the domain of the plate with boundary $S_{D} $.

It is known [1] that the function $\omega \,\left( {x_{1} x_{2} ,\,t}
\right)$ describing across vibrations of the plate satisfies equation

\begin{equation}
\label{eq1}
\omega _{x_{1} x_{1} x_{1} x_{1}}  + 2\omega _{x_{1} x_{1} x_{2} x_{2}}  +
\omega _{x_{2} x_{2} x_{2} x_{2}}  + \omega _{tt} = 0 \quad .
\end{equation}

Assuming the process stabilized the solution - eigen-vibration is seeking as

\[
\omega (x_1 ,x_2 ,t) = u(x_1 ,x_2 )\cos \lambda t,
\]

\noindent
where $\lambda $-is an eigen-frequency.

Substituting to (\ref{eq1}) we get

\begin{equation}
\label{eq2}
\Delta ^{2} = \lambda u,\,\,\,\,\,x \in D,\,\,\,
\end{equation}

\noindent where $\Delta ^{2} = \Delta \Delta $, $ \Delta $-
Laplace operator.

For different cases different boundary conditions may be put. The object
under investigation is the freezed plate with boundary conditions

\begin{equation}
\label{eq3}
u = 0,\,\,\,\,\frac{{\partial u}}{{\partial n}} = 0,\,\,\,\,\,x \in S_{D} ,
\end{equation}

\noindent here $D \in R^{2}$ -bounded convex domain with boundary
$S_{D} \in C^{2}$. It is known [5], that eigen-frequency
$\lambda_j $ - is positive and may be numbered as $0 \le \lambda
_{1} \le \lambda _{2} \le ....$ The set of all convex bounded
domains $D \in E^{2}$ we denote by $M$. Let

\[
K = \left\{ {D \in M:S_{D} \in \dot {C}^{2}} \right\},
\]

\noindent where $\dot C^2 $ is a class of the piece-wise twice
continuous differentiable functions.

 $K$ may be defined by various ways, for example, by fixing of area or length of
boundary, or by condition $D_{1} \subset D \subset D_{2} $, where $D_{1}
,\,D_{2} \in M$ are given domains. The problem is: To find a domain $D \in
K$ , such that

\begin{equation}
\label{eq4} \frac{{\left| {\Delta u_j (x)} \right|^2 }}{{\lambda
_j^{} }} = s_j (x),x \in S_D ,j = 1,2,... ,
\end{equation}

\noindent
where $u_{j} \left( {x} \right)$ and $\lambda _{j} $ are eigen-vibration and
eigen-frequency of the problem (\ref{eq2})-(\ref{eq3}) in the domain D correspondingly,
$s_{j} \left( {x} \right)$-given functions. Note that $s_{j} \left( {x}
\right),\,j = 1,2,...$ we call s-functions of the problem (\ref{eq2}),(\ref{eq3}).

In [6] the following formula is obtained for the eigenfrequency of the
freezed plate under across vibrations

\begin{equation}
\label{eq5} \lambda _j  = \frac{1}{4}\mathop {\max }\limits_{u_j }
\int\limits_{S_D } {\left| {\Delta u_j (\xi )} \right|^2 P_D
(n(\xi ))ds,}
\end{equation}

\noindent where $P_{D} \left( {x} \right) = \mathop
{max}\limits_{l \in D} \left( {l,x} \right),\,\,\,\,x \in E^{n} -
$is a support function of $D$, and $max$ is taken over all
eigen-vibrations $u_{j} $ corresponding to eigen-frequency
$\lambda _j $ in the case of its multiplicity. As we see, the
boundary values of the function $\left| {\Delta u_{j} \left( {x}
\right)} \right|^{2}$ unequivocally define $\lambda _{j} $. From
(\ref{eq5}) considering (\ref{eq4}) we get

\begin{equation}
\label{eq6}
\int\limits_{S_{D}}  {S_{j}}  \left( {\xi}  \right)P_{D} \left( {n\left( {x}
\right)} \right)ds = 4,\,\,\,j = 1,2,...
\end{equation}

This relation is basic for solving of the considering problem by given
$s$-functions.

Now we prove the lemma that will be used in further.

\textbf{Lemma1.} Let $f(x)$ be continuous function on $S_{B} $.
Then for any

\[
D_{1} ,\,\,D_{2} \in K
\]

\begin{equation}
\label{eq7} \int\limits_{S_{D_1  + D_2 } } {f(n} (\xi ))d\xi  =
\int\limits_{S_{D_1 } } {f(n(\xi ))d\xi  + } \int\limits_{S_{D_2 }
} {f(n(\xi ))d\xi } ,
\end{equation}

\noindent
where $D_{1} + D_{2} $ is taken in the sence of Minkovsky i.e.
$
D_1  + D{}_2^{}  = \left\{ {x:x = x_1  + x_2 ,x{}_1^{}  \in D_1
,x_2  \in D_2 } \right\}
$, B-unit sphere.

\textbf{Proof.} It is known, that\textbf{}  $f\left( {x} \right)$ may be
continuously, positively defined extended over all the space and presented
as a limit of the difference of two convex functions [7]

\begin{equation}
\label{eq8}
f\left( {x} \right) = \mathop {lim}\limits_{n \to \infty}  \left[ {g_{n}
\left( {x} \right) - h_{n} \left( {x} \right)} \right].
\end{equation}

Not disrupting integrity we can write

\begin{equation}
\label{eq9}
f\left( {x} \right) = g\left( {x} \right) - h\left( {x} \right),
\end{equation}

\noindent
where $g\left( {x} \right),\,h\left( {x} \right)$ are convex, positively
defined functions.

It is known [8] that for any continuous, convex,
positively-defined function $ P(x) $ there exists the only convex
bounded domain $D$ such, that $P\left( {x} \right)$ is a support
function of $D$, i.e. $P\left( {x} \right) = P_{D} \left( {x}
\right)$. The opposed statement is also true. It is also known
that $D$ is a subdifferential of its support function at the point
$x = 0$

\[
D = \partial P\left( {0} \right) = \left\{ {l \in E^{n}:P\left( {x} \right)
\ge \left( {l,x} \right),\,\forall x \in E^{n}} \right\}.
\]

So there exist the domains $G$and $H$ such that

\begin{equation}
\label{eq10} g(x) = P_G (x),h(x) = P_H (x),x \in B
\end{equation}

Considering (\ref{eq9}), (\ref{eq10}) we get

\begin{equation}
\label{eq11}
\begin{array}{l}
 \int\limits_{S_{D_1  + D_2 } } {f(n(x))ds = } \int\limits_{S_{D_1  + D_2 } }
  {\left[ {g(n(x)) - h(n(x)))ds} \right] = }  \\
  = \int\limits_{S_{D_1  + D_2 } } {P_G (n(x))ds - }
  \int\limits_{S_{D_1  + D_2 } } {P_H (n(x))ds.}  \\
 \end{array}
\end{equation}

As [6] for any $D_{1} $,$D_{2} $

\begin{equation}
\label{eq12} \int\limits_{S_{D_1 } } {P_{D_2 } (n(x))ds = }
\int\limits_{S_{D_2 } } {P_{D_1 } (n(x)))ds}
\end{equation}

\noindent from (\ref{eq11}) one may obtain $ \int\limits_{S_{D_1 +
D_2 } } {f(n(x))ds = } \int\limits_{S_G } {P_{D_1  + D_2 }
(n(x))ds}  - \int\limits_{S_H } {P_{D_1  + D_2 } (n(x))ds} $.

As $P_{D_{1} + D_{2}}  \left( {x} \right) = P_{D_{1}}  \left( {x} \right) +
P_{D_{2}}  \left( {x} \right)$ [8] applying (\ref{eq12}) again we get (\ref{eq7}).

Lemma is proved.

Now we investigate the main problem of the work- reconstraction of $D$ by
given $s$- functions.

Let $B \subset E^{2}$ and $ S_B $- its boundary. By $\varphi _{k}
\left( {x} \right),\,\,k = 1,2,..$ we denote some basis in
$C\left( {S_{B}}  \right)$-space of continuous functions on $S_{B}
$. These functions may be continuously, positive homogeneously
extended to Â. It may be done as:

\begin{equation}
\label{eq13} \tilde \varphi _k (x) = \left\{ \begin{array}{l}
 \varphi _k \left( {\frac{x}{{\left\| x \right\|}}} \right) \cdot \left\| x \right\|,x \in B,x \ne 0, \\
 0,x = 0. \\
 \end{array} \right.
\end{equation}

One may test, that these functions are continuous and satisfy the positive
homogeneity condition

\[
\tilde \varphi _j (\alpha x) = \alpha \tilde \varphi _k (x),\alpha
> 0
\]

Not disrupting integrity we can denote $\tilde {\varphi} _{j} \left( {x}
\right)$ by $\varphi _{k} \left( {x} \right)$.

Thus we obtain the system of continuous, positive homogeneous functions
defined on $B$.

As we noted above each positive homogeneous, continuous function $\varphi
_{j} \left( {x} \right)$may be presented in the form

\begin{equation}
\label{eq14} \varphi _k (x) = \mathop {\lim }\limits_{n \to \infty
} \left[ {g_n^k (x) - h_n^k (x)} \right]
\end{equation}

\noindent
and there exist satisfying above mentioned properties domains $G_{n}^{k} $
and $H_{n}^{k} $ such, that

\begin{equation}
\label{eq15} g_{n}^{k} \left( {x} \right) = P_{G_{n}^{k}}  \left(
{x} \right), \quad  h_n^k (x) = P_{H_n^k } (x).
\end{equation}

These domains we call basic domains. Substituting these into (\ref{eq14}) we get

\begin{equation}
\label{eq16}
\varphi _{k} \left( {x} \right) = \mathop {lim}\limits_{n \to \infty}
\left[ {P_{G_{n}^{k}} ^{k} \left( {x} \right) - P_{H_{n}^{k}} ^{k} \left(
{x} \right)} \right].
\end{equation}

For similarity let's suppose that

\begin{equation}
\label{eq17}
\varphi _{k} \left( {x} \right) = P_{G^{k}}^{} \left( {x} \right) -
P_{H^{k}} \left( {x} \right),
\end{equation}

\noindent
where $G^{k}$and $H^{k}$ are closed, bounded convex domains.

As $n(x) \in S_B$, for any $x \in S_{D} $, we can decomposite
$P_{D} \left( {x} \right)$, $x \in S_{B} $ by basic functions
$\varphi _{k} \left( {x} \right)$

\begin{equation}
\label{eq18} P_D (x) = \sum\limits_{k = 1}^\infty  {\alpha _k
\varphi _k (x),x \in S_B } ,\alpha  \in R .
\end{equation}

Considering (\ref{eq17}) from this one may get

\begin{equation}
\label{eq19}
P_{D} \left( {x} \right) = \sum\limits_{k = 1}^{\infty}  {\alpha _{k} \left(
{P_{G^{k}} \left( {x} \right) - P_{H^{k}} \left( {x} \right)} \right)}
,\,\,\,x \in S_{B} .
\end{equation}

The set of all indexes for which $\alpha _{k} \ge 0\left( {\alpha _{k} <
0} \right)$ denote by $I^{ + \,\,}\left( {I^{ -} } \right)$.

Then the relation (\ref{eq19}) may be written as

\begin{equation}
\label{eq20} P_D (x) = \sum\limits_{k = 1}^\infty  {\alpha _k
\varphi _k (x),x \in S_B } ,\alpha  \in R
\end{equation}

From last taking into account the properties of support functions [8] we
obtain

\[
D - \sum\limits_{k \in I^{ -} }^{} {\alpha _{k} G^{k} + \sum\limits_{k \in
I^{ +} } {\alpha _{k} H^{k} =} }  \sum\limits_{k \in I^{ +} }^{} {\alpha
_{k} G^{k} - \sum\limits_{k \in I^{ -} } {\alpha _{k} H^{k}}}  .
\]

The use (\ref{eq20}) and the lemma give

\[
\begin{array}{l}
 \int\limits_{S_D } {s_j } (\xi )P_D (n(\xi )d\xi  +
 \int\limits_{\sum\limits_{k \in I^ -  } {( - \alpha _k )S_{G^k } } }
 {s_j (\xi )P_D (n(\xi ))d\xi  + }  \\
  + \int\limits_{\sum\limits_{k \in I^ +  } {\alpha _k } S_{H^k } } {s_j }
  (\xi )P_D (n(\xi )d\xi  = \int\limits_{\sum\limits_{k \in I^ +  } {\alpha _k S_{G^k } } }
  {s_j (\xi )P_D (n(\xi ))d\xi }  +  \\
  + \int\limits_{\sum\limits_{k \in I^ -  }
  {( - \alpha _k )S_{H^k } } } {s_j } (\xi ))P_D (n(\xi ))d\xi . \\
 \end{array}
\]

From this considering (\ref{eq7}) we have

\begin{equation}
\label{eq21}
\begin{array}{l}
 \int\limits_{S_D } {s_j } (\xi )P_D (n(\xi )d\xi  = \sum\limits_{k = 1}^\infty  {\alpha _k \left[ {\int\limits_{S_{G^k } } {s_j } } \right.} (\xi )P_D (n(\xi ))d\xi  -  \\
 \left. { - \int\limits_{S_{H^k } } {s_j } (\xi )P_D (n(\xi )d\xi } \right] = 4. \\
 \end{array}
\end{equation}

Substituting here (\ref{eq18}) one may get

\begin{equation}
\label{eq22}
\sum\limits_{k,m = 1}^{\infty}  {A_{k,m} \left( {j} \right)\alpha _{k}
\alpha _{m} = 4,\,\,\,\,\,j = 1,2,...} ,
\end{equation}

\noindent
where

\[
\begin{array}{l}
 A_{k,m} (j) = \int\limits_{S_{G^k } } {s_j (x)\left[ {P_{G^m } (n(x)) - P_{H^m } (n(x))} \right]ds}  -  \\
  - \int\limits_{S_{H^k } } {s_j (x)\left[ {P_{G^m } (n(x)) - P_{H^m } (n(x))} \right]ds} . \\
 \end{array}
\]

The equation (\ref{eq22}) has, generally speaking, no only solution. The function
$P_{D} \left( {x} \right)$ is reconstructed by the help of the solution of
this equation using (\ref{eq18}).

As we noted above domain $D$ is unequivocally defined by its support
function $P_{D} \left( {x} \right)$. Suppose that (\ref{eq22}) has the only solution
providing convexity of the support function of $D$.

Let's show that the expressions $\frac{{\left| {\Delta u_j (\xi )}
\right|^2 }}{{\lambda _j }},j = 1,2,...,$ for the problem
(\ref{eq2}), (\ref{eq3}) in the reconstructed by the help of this
solution using (\ref{eq18}) domain $D$ indeed are $s$-functions.
Really, if $D$ is a domain in which the problem (\ref{eq2}),
(\ref{eq3}) has given by formulae (\ref{eq18}) $s$-functions then
decomposition $\overline D $ by formulae (\ref{eq18}) and making
above done transformations we get the equation (\ref{eq22}) with
the same coefficients.

From the assumption that this equation has the only solution, it follows
$\overline {D} = D$.

If (\ref{eq22}) has more than one solution then the searching domain is in among the
ones, constracted by (\ref{eq18}) using these solutions, providing convexity of
$P\left( {x} \right)$.\textbf{}

\textbf{References}

1. S.H. Gould. Variational Methods for Eigenvalue Problems. Oxford
Univerrsty Press, 1966, 328 p.

2. S.A. Avdonin, M.I. Belishev. Boundary control and inverse problem for
non-selfadjoint Sturm-Liouville operator (BC-method), Control and
Cybernetics, 25, 1996, p.429-440.

3. Pesaint Fabrisio, Zolesio Jean-Paul. Derivees par rapport au domaine des
valeurs propres du Laplacien. C.r.Acad. sci.
Ser,1.-1995-321,¹10.,p.1337-1310.

4. J.Elschner, G.Schimdt, M.Yamamoto. An inverse problem in periodic
diffractive optics: global uniqueness with a single wave number. Inverse
problems, 2003, 19, p.779-787.

5. Vladimirov V.S. Equations of mathematical physics. M.: Nauka\textbf{,}
1988, 512 p.

6. A.A. Niftieyev, Y.S. Gasimov. Control by boundaries and eigenvalue
problems with variable domain. Publishing house ``BSU'', 2004, 185 p.

7. D.M. Burago, V.A.Zalgarmer. Geometrical inequalities. M.: Nauka, 1981,
400 p.

8. V.F. Demyanov, A.M. Rubinov. Basises of non-smooth analyses and
quasidifferential calculas. M.: Nauka, 1990.

9. Gasimov Y.S. Niftiev A.A. On a minimization of the eigenvalues of
Shrodinger operator over domains. Doclady RAS\textbf{,} 2001, v. 380, ¹3, p.
305-307.\textbf{}

10. Y.S. Gasimov. On some properties of the eigenvalues by the variation of
the domain. Mathematical physics, analyses, geometry, 2003, v.10, ¹2,
p.249-255.

\end{document}